\documentclass[11pt]{article}

\usepackage{amsmath,amssymb,amsfonts}
\usepackage{pstricks,pst-node}

\newtheorem{thm}{Theorem}
\newtheorem{prp}[thm]{Proposition}
\newtheorem{lmm}[thm]{Lemma}   
\newtheorem{crl}[thm]{Corollary}

\def\sf#1{\textsf{#1}}
\def\e_ref#1{(\ref{#1})}

\def\ov#1{\overline{#1}}
\def\ti#1{\tilde{#1}}

\def\BE#1{\begin{equation}\label{#1}}
\def\EE{\end{equation}}
\def\lr#1{\langle{#1}\rangle}

\def\lra{\longrightarrow}

\def\eset{\emptyset}
\def\i{\infty}

\def\be{\beta}

\def\ep{\epsilon}
\def\io{\iota}

\def\na{\nabla}

\def\si{\sigma}

\def\bI{\mathbb I}
\def\cJ{\mathcal J}
\def\R{\mathbb R}

\def\Z{\mathbb Z}

\def\De{\Delta}

\def\td{\textnormal d}
\def\te{\textnormal e}

\def\Bd{\textnormal{Bd}~\!}

\def\id{\textnormal{id}}

\def\Int{\textnormal{Int}~\!}

\def\St{\textnormal{St}}
\def\sd{\textnormal{sd}\,}

\begin{document}

\title{On Transverse Triangulations}
\author{Aleksey Zinger\thanks{Partially supported by DMS grant 0846978}}
\date{\today}
\maketitle

\begin{abstract}
We show that every smooth manifold admits a smooth triangulation transverse to a given smooth map.
This removes the properness assumption on the smooth map used
in an essential way in Scharlemann's construction~\cite{Sc}.
\end{abstract}

\section{Introduction}

\noindent
For $l\!\in\!\Z^{\ge0}$, let $\De^l\!\subset\!\R^l$ denote the standard $l$-simplex.
If $|K|\!\subset\!\R^N$ is a geometric realization of 
a simplicial complex~$K$ in the sense of \cite[Section~3]{Mu2},
for each $l$-simplex $\si$ of $K$ there is an injective 
linear map $\io_{\si}\!:\De^l\!\lra\!|K|$ taking $\De^l$ to~$|\si|$.\footnote{i.e.~$\io_{\si}$
takes the vertices of $\De^l$ to the vertices of $|\si|$ and is linear between them,
as in \cite[Footnote~5]{Z}}
If $X$ is a smooth manifold, a topological embedding $\mu\!:\De^l\!\lra\!X$ 
is a \sf{smooth embedding} if there exist an open neighborhood $\De_{\mu}^l$ of $\De^l$
in $\R^l$ and a smooth embedding $\ti\mu\!:\De_{\mu}^l\!\lra\!X$ so that $\ti\mu|_{\De^l}\!=\!\mu$.
A \sf{triangulation} of a smooth manifold~$X$ is a pair
$T\!=\!(K,\eta)$ consisting of a simplicial complex
and a homeomorphism $\eta\!:|K|\!\lra\!X$ such that 
$$\eta\!\circ\!\io_{\si}\!:\De^l\lra X$$
is a smooth embedding for every $l$-simplex $\si$ in $K$ and
$l\!\in\!\Z^{\ge0}$.
If $T\!=\!(K,\eta)$ is a triangulation of~$X$ and $\psi\!:X\!\lra\!X$ is a diffeomorphism,
then $\psi_*T\!=\!(K,\psi\!\circ\!\eta)$ is also a triangulation of~$X$.

\begin{thm}\label{main_thm}
If $X,Y$ are smooth manifolds and $h\!:Y\!\lra\!X$ is a smooth map,
there exists a triangulation $(K,\eta)$ of $X$ 
such that $h$ is transverse to $\eta|_{\Int\si}$
for every simplex $\si\!\in\!K$.
\end{thm}

\noindent
This theorem is stated in \cite{Z} as Lemma 2.3 and described as an obvious fact. 
As pointed out to the author by Matthias Kreck, Scharlemann~\cite{Sc} proves Theorem~\ref{main_thm}
under the assumption that the smooth map~$h$ is proper, and his argument makes use
of this assumption in an essential way.
For the purposes of~\cite{Z}, a transverse $C^1$-triangulation would suffice, and 
the existence of a such triangulation is fairly evident from
the point of view of Sard-Smale Theorem \cite[(1.3)]{Sm}.
In this note we give a detailed proof of Theorem~\ref{main_thm} as stated above,
using Sard's theorem \cite[Section~2]{Mi97}.\\

\noindent
The author would like to thank M.~Kreck for detailed comments and suggestions on \cite{Z}
and earlier versions of this note, as well as D.~McDuff and J.~Milnor for related discussions.

\section{Outline of proof of Theorem~\ref{main_thm}}

\noindent
If $K$ is a simplicial complex, we denote by $\sd K$ the barycentric subdivision of~$K$.
For any nonnegative integer~$l$, let $K_l$ be  \sf{the $l$-th skeleton of~$K$},
i.e.~the subcomplex of~$K$ consisting of the simplices in $K$ of dimension at most~$l$. 
If $\si$ is a simplex in a simplicial complex $K$ with geometric realization~$|K|$, let
$$\St(\si,K)=\bigcup_{\si\subset\si'}\Int\si'$$
be \sf{the star of $\si$ in $K$}, as in \cite[Section~62]{Mu2}, 
and $b_{\si}\in\sd K$ the barycenter of~$\si$.
The main step in the proof of Lemma~2.3 is the following observation.

\begin{prp}\label{main_prp}
Let $h\!:Y\!\lra\!X$ be a smooth map between smooth manifolds.  
If $(K,\eta)$ is a triangulation of~$X$
and $\si$ is  an $l$-simplex in~$K$, there exists a diffeomorphism
$\psi_{\si}\!:X\!\lra\!X$ restricting to the identity outside of $\eta(\St(b_{\si},\sd K))$
so that $\psi_{\si}\!\circ\!\eta|_{\Int\si}$ is transverse to~$h$.
\end{prp}

\noindent
If $\si$ and $\si'$ are two distinct simplices in $K$ of the same dimension $l$,
\BE{Stinter_e}\St(b_{\si},\sd K)\cap\St(b_{\si'},\sd K)=\eset.\EE
Since $\psi_{\si}$ is the identity outside~of $\eta(\St(b_{\si},\sd K))$ and
the collection $\{\St(b_{\si},\sd K)\}$ is locally finite, 
the composition $\psi_l\!:X\!\lra\!X$ of all diffeomorphisms $\psi_{\si}\!:X\!\lra\!X$ 
taken over all $l$-simplices~$\si$ in~$K$ is a 
well-defined diffeomorphism of~$X$.\footnote{The locally finite property implies
that the composition of these diffeomorphisms in any order is a diffeomorphism;
by \e_ref{Stinter_e}, these diffeomorphisms commute and so
the composition is independent of the order.}
Since $\psi_l\!\circ\!\eta|_{|\si|}=\psi_{\si}\!\circ\!\eta|_{|\si|}$
for every $l$-simplex $\si$ in~$K$,
we obtain the following conclusion from Proposition~\ref{main_prp}. 

\begin{crl}\label{main_crl}
Let $h\!:Y\!\lra\!X$ be a smooth map between smooth manifolds.  
If $(K,\eta)$ is a triangulation of~$X$,
for every $l\!=\!0,1,\ldots,\dim X$, there exists a diffeomorphism
$\psi_l\!:X\!\lra\!X$ restricting to the identity on $\eta(|K_{l-1}|)$
so that $\psi_l\!\circ\!\eta|_{\Int\si}$ is transverse to~$h$
for every $l$-simplex~$\si$ in~$K$.
\end{crl}

\noindent
This corollary implies Theorem~\ref{main_thm}.
By~\cite[Chapter~II]{Mu1}, $X$ admits a triangulation $(K,\eta_{-1})$.
By induction and Corollary~\ref{main_crl}, for each $l=0,1,\ldots,\dim X-1$
there exists a triangulation $(K,\eta_l)\!=\!(K,\psi_l\!\circ\!\eta_{l-1})$ of~$X$ 
which is transverse 
to $h$ on every open simplex in $K$ of dimension at most~$l$.

\section{Proof of Proposition~\ref{main_prp}}

\begin{lmm}\label{rho_lmm}
For every $l\in\Z^+$, there exists a smooth function $\rho_l\!:\R^l\lra\bar\R^+$ 
such that 
$$\rho_l^{-1}(\R^+)=\Int\De^l.$$
\end{lmm}

\noindent
{\it Proof:} Let $\rho\!:\R\lra\R$ be the smooth function given by 
$$\rho(r)=\begin{cases} \te^{-1/r},&\hbox{if}~r>0;\\
0,&\hbox{if}~r\le0. \end{cases}$$
The smooth function $\rho_l\!:\R^l\!\lra\!\R$ given by
$$\rho_l(t_1,\ldots,t_n)=\rho\bigg(1-\sum_{i=1}^{i=l}t_i\bigg)\cdot\prod_{i=1}^{i=l}\rho(t_i)$$
then has the desired property.

\begin{lmm}\label{slope_lmm}
Let $(K,\eta)$ be a triangulation of a smooth manifold~$X$ and
$\si$ an $l$-simplex in~$K$.
If
$$\ti\mu_{\si}\!: \De_{\si}^l\times\R^{m-l}\lra U_{\si}\subset X$$
is a diffeomorphism onto an open neighborhood $U_{\si}$ of $\eta(|\si|)$ in $X$
such that $\ti\mu_{\si}(t,0)=\eta(\io_{\si}(t))$ for all $t\!\in\!\De_{\si}$,
there exists $c_{\si}\!\in\!\R^+$ such that 
$$\big\{(t,v)\!\in\!(\Int\De^l)\!\times\!\R^{m-l}\!:\,|v|\!\le\!c_{\si}\rho_l(t)\big\}
\subset\ti\mu_{\si}^{-1}\big(\eta(\St(b_{\si},\sd K))\big).$$
\end{lmm}

\noindent
{\it Proof:} It is sufficient to show that there exists $c_{\si}\!>\!0$ such that 
$$\big\{(t,v)\!\in\!(\Int\De^l)\!\times\!\R^{m-l}\!:\,|v|\!\le\!c_{\si}\rho_l(t)\big\}
\subset\ti\mu_{\si}^{-1}\big(\eta(\St(\si,K))\big).
\footnote{If $K'$ is the subdivision of $K$ obtained by adding the vertices $b_{\si'}$ 
with $\si'\supsetneq\si$, then $\St(b_{\si},\sd K)=\St(\si,K')$.}$$
We assume that $0\!<\!l\!<m$.
Suppose $(t_p,v_p)\in(\Int\De^l)\!\times\!(\R^{m-l}\!-\!0)$ is a sequence such that 
\BE{setup_e0}(t_p,v_p)\not\in  \ti\mu_{\si}^{-1}\big(\eta(\St(\si,K))\big),
\qquad |v_p|\le \frac{1}{p}\rho_l(t_p).\EE
Since $\eta(\St(\si,K))$ is an open neighborhood of $\eta(\Int\si)$ in $X$, 
by shrinking $v_p$ and passing to a subsequence we can assume that 
\BE{setup_e}
(t_p,v_p)\in\ti\mu_{\si}^{-1}\big(\eta(|\tau'|)\big)\subset\ti\mu_{\si}^{-1}\big(\eta(|\tau|)\big)\EE
for an $m$-simplex $\tau$ in $K$ and a face $\tau'$ of $\tau$ so that $\si\!\not\subset\!\tau'$,  
$\tau'\!\not\subset\!\si$, and $\si\subset\tau$.
Let $\io_{\tau}\!: \De^m\!\lra\!|K|$ be
an injective linear map taking $\De^m$ to~$|\tau|$ so that 
\BE{tausi_e}\io_{\tau}^{-1}(|\si|)=\De^m\cap\R^l\!\times\!0\subset\R^l\!\times\!\R^{m-l}\,,
\qquad \io_{\tau}^{-1}(|\tau'|)=\De^m\cap0\!\times\!\R^{m-1}\subset\R^1\!\times\!\R^{m-1}\,.\EE
Choose a smooth embedding $\mu_{\tau}\!:\De_{\tau}^m\!\lra\!X$  from an open neighborhood of 
$\De^m$ in~$\R^m$ such that \hbox{$\mu_{\tau}|_{\De^m}\!=\!\eta\circ\!\io_{\tau}$}.
Let $\phi$ be the first component of the diffeomorphism
$$\mu_{\tau}^{-1}\!\circ\!\ti\mu_{\si}\!: \ti\mu_{\si}^{-1}\big(\mu_{\tau}(\De^m_{\tau})\big)
\lra\mu_{\tau}^{-1}\big(\mu_{\si}(\De^l_{\si}\!\times\!\R^{m-l})\big)
\subset\R^1\!\times\!\R^{m-1}\,.$$
By \e_ref{setup_e}, the second assumption in~\e_ref{tausi_e}, the continuity of $\td\phi$,
and the compactness of~$\De^l$,  
\BE{phit_e}\big|\phi(t_p,0)\big|=\big|\phi(t_p,0)-\phi(t_p,v_p)\big| \le C|v_p|
\qquad\forall\,p,\EE
for some $C\!>\!0$.
On the other hand, by the first assumption in \e_ref{tausi_e},
the vanishing of~$\rho_l$ on~$\Bd\De^l$, the continuity of~$\td\rho_l$, 
and the compactness of~$\De^l$, 
\BE{rhol_e} \big|\rho_l(t_p)\big|\le C\big|\phi(t_p,0)\big| \qquad\forall\,p,\EE
for some $C\!>\!0$.
The second assumption in~\e_ref{setup_e0}, \e_ref{phit_e}, and~\e_ref{rhol_e} 
give a contradiction for $p>C^2$.

\begin{lmm}\label{SS_lmm}
Let $h\!:Y\!\lra\!X$ be a smooth map between smooth manifolds,
$(K,\eta)$ a triangulation of~$X$, $\si$ an $l$-simplex in~$K$, and
$$\ti\mu_{\si}\!: \De_{\si}^l\times\R^{m-l}\lra U_{\si}\subset X$$
a diffeomorphism onto an open neighborhood $U_{\si}$ of $\eta(|\si|)$ in $X$
such that $\ti\mu_{\si}(t,0)=\eta(\io_{\si}(t))$ for all $t\!\in\!\De_{\si}$.
For every $\ep\!>\!0$, there exists
$s_{\si}\!\in\!C^{\i}(\Int\De^l;\R^{m-l})$ so that the map
\BE{deform_e} \ti\mu_{\si}\circ(\id,s_{\si})\!:\Int\De^l\lra X\EE
is transverse to $h$,
\BE{SS_e2}
\big|s_{\si}(t)\big|<\ep^2\rho_l(t)~~~
\forall\,t\!\in\!\Int\De^l\,,\qquad
\lim_{t\lra\Bd\De^l}\rho_l(t)^{-i}\big|\na^js_{\si}(t)\big|=0~~\forall i,j\in\Z^{\ge0},  \EE
where $\na^js_{\si}$ is the multi-linear functional determined by the $j$-th
derivatives of $s_{\si}$.
\end{lmm}

\noindent
{\it Proof:} The smooth map
$$\phi\!:\Int\De^l\times\R^{m-l}\lra X, \qquad 
\phi(t,v)=\ti\mu_{\si}\big(t,\te^{-1/\rho_l(t)}v\big),$$
is a diffeomorphism onto an open neighborhood $U_{\si}'$ of $\eta(\Int\si)$ in $X$.
The smooth map~\e_ref{deform_e} with $s_{\si}\!=\!\te^{-1/\rho_l(t)}v$ is transverse to~$h$ 
if and only~if $v\!\in\!\R^{m-l}$ is a regular value of the smooth map
$$\pi_2\!\circ\!\phi^{-1}\!\circ\!h\!: h^{-1}(U_{\si}')\lra\R^{m-l}\,,$$
where $\pi_2\!:\Int\De^l\!\times\!\R^{m-l}\!\lra\!\R^{m-l}$ is the projection 
onto the second component.
By Sard's Theorem, the set of such regular values 
is dense in~$\R^{m-l}$.
Thus, the map~\e_ref{deform_e} with $s_{\si}\!=\!\te^{-1/\rho_l(t)}v$  is transverse to~$h$ 
for some $v\!\in\!\R^{m-l}$ with $|v|\!<\!\ep^2$.
The second statement in~\e_ref{SS_e2} follows from $\rho_l|_{\Bd\De^l}\!=\!0$.

\begin{crl}\label{ext_crl}
Let $h\!:Y\!\lra\!X$ be a smooth map between smooth manifolds,
$(K,\eta)$ a triangulation of~$X$, $\si$ an $l$-simplex in~$K$, and
$$\ti\mu_{\si}\!: \De_{\si}^l\times\R^{m-l}\lra U_{\si}\subset X$$
a diffeomorphism onto an open neighborhood $U_{\si}$ of $\eta(|\si|)$ in $X$
such that $\ti\mu_{\si}(t,0)=\eta(\io_{\si}(t))$ for all $t\!\in\!\De_{\si}$.
For every $\ep\!>\!0$, there exists  a diffeomorphism~$\psi_{\si}'$ of $\De^l_{\si}\!\times\!\R^{m-l}$
restricting to the identity outside~of 
$$\big\{(t,v)\!\in\!(\Int\De^l)\!\times\!\R^{m-l}\!:\,|v|\!\le\!\ep\rho_l(t)\big\}$$
so that the map $\ti\mu_{\si}\!\circ\psi_{\si}'|_{\Int\De^l\times0}$ is transverse to~$h$.
\end{crl}

\noindent
{\it Proof:} Choose $\be\!\in\!C^{\i}(\R;[0,1])$ so that 
$$\be(r)=\begin{cases} 1,&\hbox{if}~r\le\frac{1}{2};\\ 
0,&\hbox{if}~r\ge1.\end{cases}$$
Let $C_{\be}\!=\!\sup_{r\in\R}|\be'(r)|$.
With $s_{\si}$ as provided by Lemma~\ref{SS_lmm}, define 
\begin{gather*}
\psi_{\si}'\!:\De_{\si}^l\times\R^{m-l}\lra \De_{\si}^l\times\R^{m-l}  \qquad\hbox{by}\\
\psi_{\si}'(t,v)=\begin{cases}
\left(t,v+\be\left(\frac{|v|}{\ep\rho_l(t)}\right)s_{\si}(t)\right),&\hbox{if}~t\in\Int\De^l;\\
(t,v),&\hbox{if}~t\not\in\Int\De^l.  \end{cases}\end{gather*}
The restriction of this map to $(\Int\De^l)\!\times\!\R^{m-l}$ is smooth and its Jacobian is
\BE{Jac_e}\cJ\psi_{\si}'\big|_{(t,v)}=
\left(\begin{array}{cc} \bI_l& 0\\ 
\be\left(\frac{|v|}{\ep\rho_l(t)}\right)\na s_{\si}(t)
-\be'\left(\frac{|v|}{\ep\rho_l(t)}\right)\frac{|v|}{\ep\rho_l(t)}\frac{s_{\si}(t)}{\rho_l(t)}\na\rho_l
& \bI_{m-l}+\be'\left(\frac{|v|}{\ep\rho_l(t)}\right)\frac{s_{\si}(t)}{\ep\rho_l(t)}\frac{v^{tr}}{|v|}
\end{array}\right).\EE
By the first property in \e_ref{SS_e2}, this matrix is non-singular if $\ep\!<\!1/C_{\be}$.
If $W$ is any linear subspace of~$\R^{m-l}$ containing $s_{\si}(t)$,
$$\psi_{\si}'(t\!\times\!W)\subset t\!\times\!W, \qquad
\psi_{\si}'(t,v)=(t,v)~~~\forall\,v\!\in\!W~\hbox{s.t.}~|v|\ge\ep\rho_l(t).$$
Thus, $\psi_{\si}'$ is a bijection on $t\!\times\!W$, a diffeomorphism
on $(\Int\De^l)\times\!\R^{m-l}$, and a bijection on~$\De_{\si}^l\!\times\!\R^{m-l}$.\\

\noindent
Since $\be(r)\!=\!0$ for $r\!\ge\!1$, $\psi_{\si}'(t,v)\!=\!(t,v)$ unless $t\!\in\!\Int\De^l$
and $|v|\!<\!\ep\rho_l(t)$. 
It remains to show that $\psi_{\si}'$ is smooth along
$$\ov{\big\{(t,v)\!\in\!(\Int\De^l)\!\times\!\R^{m-l}\!:\,|v|\!\le\!\ep\rho_l(t)\big\}}
-(\Int\De^l)\!\times\!\R^{m-l}
=(\Bd\De^l)\!\times\!0.$$
Since $|s_{\si}(t)|\!\lra\!0$ as $t\!\lra\!\Bd\De^l$ by the first property in~\e_ref{SS_e2}, 
$\psi_{\si}'$ is continuous at all $(t,0)\!\in\!(\Bd\De^l)\!\times\!0$.
By the first property in~\e_ref{SS_e2}, 
$\psi_{\si}'$ is also differentiable at all $(t,0)\!\in\!(\Bd\De^l)\!\times\!0$,
with the Jacobian equal to~$\bI_m$.
By~\e_ref{Jac_e} and the compactness of~$\De^l$,
$$\big|\cJ\psi_{\si}'|_{(t,v)}-\bI_m\big|
\le C\big(|\na s_{\si}(t)|+\rho(t)^{-1}|s_{\si}(t)|\big)
\quad\forall\,(t,v)\!\in\!(\Int\De^l)\!\times\!\R^{m-l}$$
for some $C\!>\!0$. 
So $\cJ\psi_{\si}'$ is continuous at $(t,0)$ by  the second statement in~\e_ref{SS_e2}, 
as well as differentiable, with the differential of $\cJ\psi_{\si}'$ at $(t,0)$ equal to~$0$.
For $i\!\ge\!2$, the $i$-th derivatives of the second component of $\psi_{\si}'$
at $(t,v)\!\in\!(\Int\De^l)\!\times\!\R^{m-l}$ are linear combinations of
the terms
$$\be^{\lr{i_1}}\left(\frac{|v|}{\ep\rho_l(t)}\right)\cdot
\bigg(\frac{|v|}{\ep\rho_l(t)}\bigg)^{i_1}
\cdot\prod_{k=1}^{k=j_1}\left(\frac{\na^{p_k}\rho_l}{\rho_l(t)}\right)\cdot\frac{v_J}{|v|^{2j_2}}
\cdot\na^{i_2}s_{\si}(t)\,,$$
where $i_1,i_2,j_1,j_2\!\in\!\Z^{\ge0}$ and $p_1,\ldots,p_{j_1}\!\in\!\Z^+$ are such that 
$$i_1+(p_1+p_2+\ldots+p_{j_1}-j_1)+i_2=i, \qquad
j_1\!+\!j_2\!\le\!i_1\,,$$ 
and $v_J$ is a $j_2$-fold product of components of~$v$.
Such a term is nonzero only if $\ep\rho_l(t)/2\!<\!|v|\!<\!\ep\rho_l(t)$
or $i_1\!=\!0$ and $|v|\!<\!\ep\rho_l(t)$.
Thus, the $i$-th derivatives of $\psi_{\si}'$ at $(t,v)\!\in\!(\Int\De^l)\!\times\!\R^{m-l}$ 
are bounded~by
$$C_i\sum_{i_1+i_2\le i}\rho_l(t)^{-i_1}\big|\na^{i_2}s_{\si}(t)\big|$$
for some constant $C_i\!>\!0$.
By the second statement in~\e_ref{SS_e2}, 
the last expression approaches~$0$ as $t\!\lra\!\Bd\De^l$ 
and does so faster than $\rho_l$.
It follows that $\psi_{\si}'$ is smooth at all $(t,0)\!\in\!(\Bd\De^l)\!\times\!0$.\\

\noindent
{\it Proof of Proposition~\ref{main_prp}:}
Let $\De^l_{\si}$ be a contractible open neighborhood of $\De^l$ in $\R^l$ 
and $\mu_{\si}\!:\De^l_{\si}\!\lra\!X$ a smooth embedding so that 
$\mu_{\si}|_{\De^l}\!=\!\eta\!\circ\!\io_{\si}$.
By the Tubular Neighborhood Theorem \cite[(12.11)]{BJ}, there exist an open neighborhood~$U_{\si}$
of $\mu_{\si}(\De^l_{\si})$ in~$X$ and a diffeomorphism
$$\ti\mu_{\si}\!:\De^l_{\si}\!\times\!\R^{m-l}\lra U_{\si}
\qquad\hbox{s.t.}\qquad \ti\mu_{\si}(t,0)=\mu_{\si}(t)
~~\forall\,t\!\in\!\De^l_{\si}\,.\footnote{Since $\De^l_{\si}$ is contractible,
the normal bundle to the embedding $\mu_{\si}$ is trivial.}$$
Let $c_{\si}\!>\!0$ be as in Lemma~\ref{slope_lmm} and 
$\psi_{\si}'$ as in Corollary~\ref{ext_crl} with $\ep\!=\!c_{\si}$. 
The diffeomorphism 
$$\psi_{\si}=\ti\mu_{\si}\circ\psi_{\si}'\circ\ti\mu_{\si}^{-1}\!: U_{\si}\lra U_{\si}$$
is then the identity on $U_{\si}-\St(b_{\si},\sd K)$.
Since $\psi_{\si}$ is also the identity outside of a compact subset of~$U_{\si}$,
it extends by identity to a diffeomorphism on all of~$X$.\\

\vspace{.2in}

\noindent
{\it Department of Mathematics, SUNY, Stony Brook, NY 11790-3651\\
azinger@math.sunysb.edu}\\

\end{document}